\font\fifteenrm=cmr10 scaled\magstep2 
\font\fifteeni=cmmi10 scaled\magstep2
\font\fifteensy=cmsy10 scaled\magstep2
\font\fifteenbf=cmbx10 scaled\magstep2
\font\fifteentt=cmtt10 scaled\magstep2
\font\fifteenit=cmti10 scaled\magstep2
\font\fifteensl=cmsl10 scaled\magstep2
\font\fifteenam=msam10 scaled\magstep2
\font\fifteenbm=msbm10 scaled\magstep2
\font\fifteenex=cmex10 scaled\magstep2
\font\fifteensc=cmcsc10 scaled\magstep2 
\font\twelverm=cmr10 at 12pt
\font\twelvei=cmmi10 at 12pt
\font\twelvesy=cmsy10 at 12pt
\font\twelvebf=cmbx10 at 12pt
\font\twelvett=cmtt10 at 12pt
\font\twelveit=cmti10 at 12pt
\font\twelvesl=cmsl10 at 12pt
\font\twelveam=msam10 at 12pt
\font\twelvebm=msbm10 at 12pt
\font\twelveex=cmex10 at 12pt
\font\twelvesc=cmcsc10 at 12pt
\font\elevenrm=cmr10 scaled\magstephalf 
\font\eleveni=cmmi10 scaled\magstephalf
\font\elevensy=cmsy10 scaled\magstephalf
\font\elevenbf=cmbx10 scaled\magstephalf
\font\eleventt=cmtt10 scaled\magstephalf
\font\elevenit=cmti10 scaled\magstephalf
\font\elevensl=cmsl10 scaled\magstephalf
\font\elevenam=msam10 scaled\magstephalf
\font\elevenbm=msbm10 scaled\magstephalf
\font\elevenex=cmex10 scaled\magstephalf
\font\elevensc=cmcsc10 scaled\magstephalf
\font\tenrm=cmr10
\font\teni=cmmi10
\font\tensy=cmsy10
\font\tenbf=cmbx10
\font\tentt=cmtt10
\font\tenit=cmti10
\font\tensl=cmsl10
\font\tenam=msam10
\font\tenbm=msbm10
\font\tenex=cmex10
\font\tensc=cmcsc10
\font\ninerm=cmr9
\font\ninei=cmmi9
\font\ninesy=cmsy9
\font\ninebf=cmbx9
\font\ninett=cmtt9
\font\nineit=cmti9
\font\ninesl=cmsl9
\font\nineam=msam9
\font\ninebm=msbm9
\font\nineex=cmex9
\font\ninesc=cmcsc9
\font\eightrm=cmr8
\font\eighti=cmmi8
\font\eightsy=cmsy8
\font\eightbf=cmbx8
\font\eighttt=cmtt8
\font\eightit=cmti8
\font\eightsl=cmsl8
\font\eightam=msam8
\font\eightbm=msbm8
\font\eightex=cmex8
\font\eightsc=cmcsc8
\font\sevenrm=cmr7
\font\seveni=cmmi7
\font\sevensy=cmsy7
\font\sevenbf=cmbx7

\font\sevenam=msam7
\font\sevenbm=msbm7

\font\sixrm=cmr6
\font\sixi=cmmi6
\font\sixsy=cmsy6

\font\sixam=msam6
\font\sixbm=msbm6

\font\fiverm=cmr5
\font\fivei=cmmi5
\font\fivesy=cmsy5

\font\fiveam=msam5
\font\fivebm=msbm5

\font\fourrm=cmr5 at 4pt
\font\fouri=cmmi5 at 4pt
\font\foursy=cmsy5 at 4pt

\font\fouram=msam5 at 4pt
\font\fourbm=msbm5 at 4pt

\skewchar\twelvei='177 \skewchar\eleveni='177\skewchar\teni='177
\skewchar\ninei='177 \skewchar\eighti='177\skewchar\seveni='177 
\skewchar\sixi='177 \skewchar\fivei='177 \skewchar\fouri='177
\skewchar\twelvesy='60 \skewchar\elevensy='60 \skewchar\tensy='60
\skewchar\ninesy='60 \skewchar\eightsy='60 \skewchar\sevensy='60 
\skewchar\sixsy='60 \skewchar\fivesy='60 \skewchar\foursy='60
\newfam\itfam
\newfam\slfam
\newfam\bffam
\newfam\ttfam
\newfam\scfam
\newfam\amfam
\newfam\bmfam
\def\eightbig#1{{\hbox{$\left#1\vbox to 6.5pt{}\voidright $}}}
\def\eightBig#1{{\hbox{$\left#1\vbox to 7.5pt{}\voidright $}}}
\def\eightbigg#1{{\hbox{$\left#1\vbox to 10pt{}\voidright $}}}
\def\eightBigg#1{{\hbox{$\left#1\vbox to 13pt{}\voidright $}}}
\def\ninebig#1{{\hbox{$\left#1\vbox to 7.5pt{}\voidright $}}}
\def\nineBig#1{{\hbox{$\left#1\vbox to 8.5pt{}\voidright $}}}
\def\ninebigg#1{{\hbox{$\left#1\vbox to 11.5pt{}\voidright $}}}
\def\nineBigg#1{{\hbox{$\left#1\vbox to 14.5pt{}\voidright $}}}
\def\tenbig#1{{\hbox{$\left#1\vbox to 8.5pt{}\voidright $}}}
\def\tenBig#1{{\hbox{$\left#1\vbox to 9.5pt{}\voidright $}}}
\def\tenbigg#1{{\hbox{$\left#1\vbox to 12.5pt{}\voidright $}}}
\def\tenBigg#1{{\hbox{$\left#1\vbox to 16pt{}\voidright $}}}
\def\elevenbig#1{{\hbox{$\left#1\vbox to 9pt{}\voidright $}}}
\def\elevenBig#1{{\hbox{$\left#1\vbox to 10.5pt{}\voidright $}}}
\def\elevenbigg#1{{\hbox{$\left#1\vbox to 14pt{}\voidright $}}}
\def\elevenBigg#1{{\hbox{$\left#1\vbox to 17.5pt{}\voidright $}}}
\def\twelvebig#1{{\hbox{$\left#1\vbox to 10pt{}\voidright $}}}
\def\twelveBig#1{{\hbox{$\left#1\vbox to 11pt{}\voidright $}}}
\def\twelvebigg#1{{\hbox{$\left#1\vbox to 15pt{}\voidright $}}}
\def\twelveBigg#1{{\hbox{$\left#1\vbox to 19pt{}\voidright $}}}
\def\fifteenbig#1{{\hbox{$\left#1\vbox to 12pt{}\voidright $}}}
\def\fifteenBig#1{{\hbox{$\left#1\vbox to 13.5pt{}\voidright $}}}
\def\fifteenbigg#1{{\hbox{$\left#1\vbox to 18pt{}\voidright $}}}
\def\fifteenBigg#1{{\hbox{$\left#1\vbox to 23pt{}\voidright $}}}
\def\voidright{\right.\nulldelimiterspace=0pt \mathsurround=0pt }
\def\fifteenpoint{
  \textfont0=\fifteenrm \scriptfont0=\twelverm \scriptscriptfont0=\tenrm
  \def\rm{\fam0 \fifteenrm}%
  \textfont1=\fifteeni \scriptfont1=\twelvei \scriptscriptfont1=\teni
  \textfont2=\fifteensy \scriptfont2=\twelvesy \scriptscriptfont2=\tensy
  \textfont3=\fifteenex \scriptfont3=\fifteenex \scriptscriptfont3=\fifteenex
  \def\it{\fam\itfam\fifteenit}\textfont\itfam=\fifteenit
  \def\sl{\fam\slfam\fifteensl}\textfont\slfam=\fifteensl
  \def\bf{\fam\bffam\fifteenbf}\textfont\bffam=\fifteenbf 
    \scriptfont\bffam=\twelvebf\scriptscriptfont\bffam=\tenbf
  \def\tt{\fam\ttfam\fifteentt}\textfont\ttfam=\fifteentt
  \def\sc{\fam\scfam\fifteensc}\textfont\scfam=\fifteensc
  \def\am{\fam\amfam\fifteenam}\textfont\amfam=\fifteenam
    \scriptfont\amfam=\twelveam\scriptscriptfont\amfam=\tenam
  \def\bm{\fam\bmfam\fifteenbm}\textfont\bmfam=\fifteenbm
    \scriptfont\bmfam=\twelvebm\scriptscriptfont\bmfam=\tenbm
  \baselineskip=21pt \rm
  \let\big=\fifteenbig\let\Big=\fifteenBig\let\bigg=\fifteenbigg
  \let\Bigg=\fifteenBigg}
\def\twelvepoint{
  \textfont0=\twelverm \scriptfont0=\ninerm \scriptscriptfont0=\sevenrm
  \def\rm{\fam0 \twelverm}%
  \textfont1=\twelvei \scriptfont1=\ninei \scriptscriptfont1=\seveni
  \textfont2=\twelvesy \scriptfont2=\ninesy \scriptscriptfont2=\sevensy
  \textfont3=\twelveex \scriptfont3=\twelveex \scriptscriptfont3=\twelveex
  \def\it{\fam\itfam\twelveit}\textfont\itfam=\twelveit
  \def\sl{\fam\slfam\twelvesl}\textfont\slfam=\twelvesl
  \def\bf{\fam\bffam\twelvebf}\textfont\bffam=\twelvebf 
    \scriptfont\bffam=\ninebf\scriptscriptfont\bffam=\sevenbf
  \def\tt{\fam\ttfam\twelvett}\textfont\ttfam=\twelvett
  \def\sc{\fam\scfam\twelvesc}\textfont\scfam=\twelvesc
  \def\am{\fam\amfam\twelveam}\textfont\amfam=\twelveam
    \scriptfont\amfam=\nineam\scriptscriptfont\amfam=\sevenam
  \def\bm{\fam\bmfam\twelvebm}\textfont\bmfam=\twelvebm
    \scriptfont\bmfam=\ninebm\scriptscriptfont\bmfam=\sevenbm
  \baselineskip=17.8pt \rm 
  \def\looselineskip{\baselineskip=18.5pt plus 1.8pt}%
  \def\tightlineskip{\baselineskip=16.5pt}%
  \def\verytightlineskip{\baselineskip=15pt}%
  \let\big=\twelvebig\let\Big=\twelveBig\let\bigg=\twelvebigg
  \let\Bigg=\twelveBigg  }
\def\elevenpoint{
  \textfont0=\elevenrm \scriptfont0=\ninerm \scriptscriptfont0=\sixrm
  \def\rm{\fam0 \elevenrm}%
  \textfont1=\eleveni \scriptfont1=\ninei \scriptscriptfont1=\sixi
  \textfont2=\elevensy \scriptfont2=\ninesy \scriptfont2=\sixsy 
  \textfont3=\elevenex \scriptfont3=\elevenex \scriptfont3=\elevenex
  \def\it{\fam\itfam\elevenit}\textfont\itfam=\elevenit
  \def\sl{\fam\slfam\elevensl}\textfont\slfam=\elevensl
  \def\bf{\fam\bffam\elevenbf}\textfont\bffam=\elevenbf
  \def\tt{\fam\ttfam\eleventt}\textfont\ttfam=\eleventt
  \def\sc{\fam\scfam\elevensc}\textfont\scfam=\elevensc
  \def\am{\fam\amfam\elevenam}\textfont\amfam=\elevenam
    \scriptfont\amfam=\nineam\scriptscriptfont\amfam=\sixam
  \def\bm{\fam\bmfam\elevenbm}\textfont\bmfam=\elevenbm
    \scriptfont\bmfam=\ninebm\scriptscriptfont\bmfam=\sixbm
  \baselineskip=15.1pt \rm
  \def\looselineskip{\baselineskip=16pt plus 1.5pt}%
  \def\tightlineskip{\baselineskip=14pt}%
  \def\verytightlineskip{\baselineskip=13pt}%
  \let\big=\elevenbig\let\Big=\elevenBig\let\bigg=\elevenbigg
  \let\Bigg=\elevenBigg  }
\def\tenpoint{
  \textfont0=\tenrm \scriptfont0=\eightrm \scriptscriptfont0=\fiverm
  \def\rm{\fam0 \tenrm}%
  \textfont1=\teni \scriptfont1=\eighti \scriptscriptfont1=\fivei
  \textfont2=\tensy \scriptfont2=\eightsy \scriptfont2=\fivesy 
  \textfont3=\tenex \scriptfont3=\tenex \scriptfont3=\tenex
  \def\it{\fam\itfam\tenit}\textfont\itfam=\tenit
  \def\sl{\fam\slfam\tensl}\textfont\slfam=\tensl
  \def\bf{\fam\bffam\tenbf}\textfont\bffam=\tenbf
  \def\tt{\fam\ttfam\tentt}\textfont\ttfam=\tentt
  \def\sc{\fam\scfam\tensc}\textfont\scfam=\tensc
  \def\am{\fam\amfam\tenam}\textfont\amfam=\tenam
    \scriptfont\amfam=\eightam \scriptscriptfont\amfam=\fiveam
  \def\bm{\fam\bmfam\tenbm}\textfont\bmfam=\tenbm
    \scriptfont\bmfam=\eightbm \scriptscriptfont\bmfam=\fivebm
  \baselineskip=14pt\rm
  \def\looselineskip{\baselineskip=14.8pt plus1.5pt}
  \def\tightlineskip{\baselineskip=12.6pt}%
  \def\verytightlineskip{\baselineskip=13pt}%
  \let\big=\tenbig\let\Big=\tenBig\let\bigg=\tenbigg\let\Bigg=\tenBigg  }
\def\ninepoint{
  \textfont0=\ninerm \scriptfont0=\sevenrm \scriptscriptfont0=\fourrm
  \def\rm{\fam0 \ninerm}%
  \textfont1=\ninei \scriptfont1=\seveni \scriptscriptfont1=\fouri
  \textfont2=\ninesy \scriptfont2=\sevensy \scriptfont2=\foursy 
  \textfont3=\nineex \scriptfont3=\nineex \scriptfont3=\nineex
  \def\it{\fam\itfam\nineit}\textfont\itfam=\nineit
  \def\sl{\fam\slfam\ninesl}\textfont\slfam=\ninesl
  \def\bf{\fam\bffam\ninebf}\textfont\bffam=\ninebf
  \def\tt{\fam\ttfam\ninett}\textfont\ttfam=\ninett
  \def\sc{\fam\scfam\ninesc}\textfont\scfam=\ninesc
  \def\am{\fam\amfam\nineam}\textfont\amfam=\nineam
    \scriptfont\amfam=\nineam\scriptscriptfont\amfam=\fouram
  \def\bm{\fam\bmfam\ninebm}\textfont\bmfam=\ninebm
    \scriptfont\bmfam=\ninebm\scriptscriptfont\bmfam=\fourbm
  \baselineskip=12.6pt\rm
  \def\tightlineskip{\baselineskip=11.5pt}
  \let\big=\ninebig\let\Big=\nineBig\let\bigg=\ninebigg
  \let\Bigg=\nineBigg  }
\def\eightpoint{
  \textfont0=\eightrm \scriptfont0=\fiverm \scriptscriptfont0=\fourrm
  \def\rm{\fam0 \eightrm}%
  \textfont1=\eighti \scriptfont1=\fivei \scriptscriptfont1=\fouri
  \textfont2=\eightsy \scriptfont2=\fivesy \scriptfont2=\foursy 
  \textfont3=\eightex \scriptfont3=\eightex \scriptfont3=\eightex
  \def\it{\fam\itfam\eightit}\textfont\itfam=\eightit
  \def\sl{\fam\slfam\eightsl}\textfont\slfam=\eightsl
  \def\bf{\fam\bffam\eightbf}\textfont\bffam=\eightbf
  \def\tt{\fam\ttfam\eighttt}\textfont\ttfam=\eighttt
  \def\sc{\fam\scfam\eightsc}\textfont\scfam=\eightsc
  \def\am{\fam\amfam\eightam}\textfont\amfam=\eightam
    \scriptfont\amfam=\eightam\scriptscriptfont\amfam=\fouram
  \def\bm{\fam\bmfam\eightbm}\textfont\bmfam=\eightbm
    \scriptfont\bmfam=\eightbm\scriptscriptfont\bmfam=\fourbm
  \baselineskip=11.2pt \rm
  \def\tightlineskip{\baselineskip=10.4pt}
  \let\big=\eightbig\let\Big=\eightBig\let\bigg=\eightbigg
  \let\Bigg=\eightBigg  }

\twelvepoint
\nopagenumbers
\hsize=6in\vsize=8.8in

\parskip=1pt plus 1pt

\newif\ifSpecialhead\Specialheadfalse
\newbox\specialheadbox

\def\specialhead #1\par{\Specialheadtrue\setbox\specialheadbox=\hbox{#1}}
\headline={{\ifSpecialhead\box\specialheadbox\global\Specialheadfalse\else
     \ifnum\pageno<0{\hfill\quad{\twelvebf\folio}}%
     \else\ifnum\pageno<2\hfill
     \else\hfill\twelvepoint\sc\firstmark\quad{\twelvebf\folio}\fi\fi\fi}}

\def\title#1\par{\bigskip{\def\cr{\par\center}\center\fifteenbf #1\par}\medskip}
\def\subtitle#1\par{\centerline{\fifteenrm #1}\medskip}
\def\author#1\par{\medskip{\def\cr{\par\center\twelvesc}\fifteensc\center#1\par}}
\def\center#1\par{\hfil #1\hfil\par}
\def\abstract.#1\par{\message{Abstract.}%
                    \medskip{\narrower\narrower\tenpoint\tightlineskip
                        \noindent{\bf Abstract.}#1\par}\medskip\noindent}
\def\tinyabstract.#1\par{\message{Abstract.}%
                    \medskip{\narrower\narrower\eightpoint\tightlineskip
                        \noindent{\bf Abstract.}#1\par}\medskip\noindent}
\def\bigabstract.#1\par{\message{Abstract.}%
                         \medskip{\narrower\narrower\tightlineskip
                         \noindent{\bf Abstract. }#1\par}\medskip\noindent}
\def\acknowledgement#1\par{\footnote{}{#1}}
\def\sectionskip{\Goodbreak\vskip 25pt plus 15pt minus 5pt}
\def\secnumber{\ifquiet
               \else\ifNoSections
                    \else\sectionsymbol\the\secno\quad\fi\fi}
\def\section#1\par{ \NoSectionsfalse\par\sectionskip\proofdepth=0\claimno=0
 \ifquiet\else\advance\secno by1\fi\toks0={#1}
 \immediate\write16{\ifquiet\else Section \the\secno\space\fi
                    \the\toks0}%
 \mark{\secnumber #1}%
 {\fifteenpoint\bf\noindent\secnumber #1}\nobreak\bigskip\quietoff
 \nobreak\noindent}
\def\quiet{\quiettrue}

\def\quietoff{\ifQUIET\else\quietfalse\fi}
\newif\ifquiet
\newif\ifQUIET
\newif\ifNoSections
\newcount\claimtype
\newcount\secno
\newcount\claimno
\newcount\subclaimno
\newcount\subsubclaimno
\newcount\subsubsubclaimno
\newcount\proofdepth
\def\subclaimnumber{\ifquiet\else\ifcase\subclaimno\or A\or B\or C\or D\or E\or
     F\or G\or H\or I\or J\or K\or L\or M\or N\or O\or P\fi\fi}
\def\subsubclaimnumber{\ifquiet\else\ifcase\subsubclaimno\or i\or ii\or iii\or 
   iv\or v\or vi\or vii\or viii\or ix\or x\or xi\or xii\or xiii\or xiv\fi\fi}
\def\subsubsubclaimnumber{\ifquiet\else\ifcase\subsubsubclaimno\or a\or b\or 
   c\or d\or e\or f\or g\or \or h\or i\or j\or k\or l\or m\or n\or o\fi\fi}
\def\claimtag{\ifquiet\else
  \ifNoSections
    \ifcase\proofdepth\the\claimno%
    \or\the\claimno.\subclaimnumber
    \or\the\claimno.\subclaimnumber.\subsubclaimnumber
    \or\the\claimno.\subclaimnumber.\subsubclaimnumber
                                                .\subsubsubclaimnumber\fi
  \else
    \ifcase\proofdepth\the\secno.\the\claimno
    \or\the\secno.\the\claimno.\subclaimnumber
    \or\the\secno.\the\claimno.\subclaimnumber.\subsubclaimnumber
    \or\the\secno.\the\claimno.\subclaimnumber.\subsubclaimnumber
                                                .\subsubsubclaimnumber\fi\fi\fi}
\secno=0\claimno=0\proofdepth=0\subclaimno=0\subsubclaimno=0\subsubsubclaimno=0
\NoSectionstrue
\newbox\qedbox
\def\claimname{\ifcase\claimtype Theorem\or Lemma\or Claim\or Corollary\or
               Question\or Definition\or Remark\or Conjecture\fi}
\def\preclaimskip{\removelastskip
    \ifcase\claimtype\goodbreak\vskip 8pt plus 4pt minus 2pt
                  \or\goodbreak\vskip 6pt plus 4pt minus 1pt
                  \or\goodbreak\vskip 5pt plus 4pt minus 1pt
                  \or\goodbreak\vskip 8pt plus 4pt minus 2pt
                  \or\vskip 7pt plus 4pt minus 2pt
                  \or\vskip 7pt plus 4pt minus 2pt
                  \or\vskip 7pt plus 4pt minus 2pt
                  \or\goodbreak\vskip 8pt plus 4pt minus 2pt\fi}
\def\postclaimskip{\ifcase\claimtype         \vskip 4pt plus 2pt minus 2pt
                                          \or\vskip 3pt plus 2pt minus 2pt
                                          \or\vskip 2pt plus 2pt minus 1pt
                                          \or\vskip 4pt plus 2pt minus 2pt
                                          \or\vskip 1pt plus 2pt 
                                          \or\vskip 4pt plus 4pt 
                                          \or\vskip 3pt plus 2pt
                                          \or\vskip 4pt plus 2pt minus 2pt\fi}
\def\claimfont{\ifcase\claimtype
                  \sl\or\sl\or\sl\or\sl\or\sl\or\rm\or\rm\or\sl\fi}
\def\advancetag{\ifcase\proofdepth\advance\claimno by1
                               \or\advance\subclaimno by1
                               \or\advance\subsubclaimno by1
                               \or\advance\subsubsubclaimno by1\fi}
\def\sayclaim#1.#2 #3\par{\ifquiet\else\advancetag\fi
    \preclaimskip\setbox1=\hbox{#1}\setbox2=\hbox{#2}%
    \toks0={#1 }
    \immediate\write16{\ifdim\wd1>0pt\the\toks0
                       \else\claimname\space\fi \claimtag.}%
    \vbox{\noindent
    {\bf\ifdim\wd1=0pt \claimname\else #1\fi\ifquiet.\else\ \claimtag{\ifNoSections.\fi}\fi}%
    \enspace{\ifdim\wd2>0pt\sc #2\enspace\fi}%
    {\claimfont #3\par}}\postclaimskip\quietoff}
\def\theorem{\claimtype=0\sayclaim}
\def\lemma{\claimtype=1\sayclaim}

\def\corollary{\claimtype=3\sayclaim}

\def\point#1. #2\par{\item{\rm #1.}#2\par}
\def\points#1\cr\par{\medskip\vbox{\let\cr=\point\point#1\par}\par}
\def\df{\it}
\def\prooffont{}
\def\proofsize{}
\def\proofindent{}
\def\proofskip{\badbreak\ifcase\claimtype    \vskip 3pt plus 2pt minus 2pt
                                          \or\vskip 2pt plus 2pt minus 2pt
                                          \or\vskip 1pt plus 2pt minus 1pt
                                          \or\vskip 3pt plus 2pt minus 2pt
                                          \or\vskip 1pt plus 2pt 
                                          \or\vskip 2pt plus 4pt 
                                          \or\vskip 1pt plus 2pt
                                          \or\vskip 3pt plus 2pt minus 2pt\fi}

\def\Goodbreak{\vskip0pt plus.5in\penalty-1000\vskip0pt plus-.5in}
\def\goodbreak{\penalty-500}
\def\badbreak{\penalty500}
\def\Badbreak{\penalty1000}
\def\proof{\message{proof}\removelastskip\Badbreak\proofskip\begingroup
  \advance\proofdepth by1
  \setbox\qedbox=\hbox{\halmos\raise2pt\hbox{\fiverm\claimname}}%
  \prooffont\proofsize\proofindent\noindent{\bf Proof: }}
\def\proofof#1:{\message{proof}\removelastskip\Badbreak\proofskip\begingroup
  \advance\proofdepth by1
  \setbox\qedbox=\hbox{\halmos\raise2pt\hbox{\fiverm#1}}%
  \prooffont\proofsize\proofindent\noindent{\bf Proof of #1: }}
\def\cite[#1]{[{\tenrm{#1}}]\message{[#1]}}
\edef\ref#1{\expandafter\global\expandafter\edef#1{\noexpand\claimtag}}
\newwrite\notes
\openout\notes=\jobname.notes
\long\def\unexpandedwrite#1#2{\def\finwrite{\write#1}%
   {\aftergroup\finwrite\aftergroup{\sanitize#2\endsanity}}}
\def\sanitize{\futurelet\next\sanswitch}
\let\stoken=\space
\def\sanswitch{\ifx\next\endsanity
  \else\ifcat\noexpand\next\stoken\aftergroup\space\let\next=\eat
   \else\ifcat\noexpand\next\bgroup\aftergroup{\let\next=\eat
    \else\ifcat\noexpand\next\egroup\aftergroup}\let\next=\eat
     \else\let\next=\copytoken\fi\fi\fi\fi \next}
\def\eat{\afterassignment\sanitize \let\next= }
\long\def\copytoken#1{\ifcat\noexpand#1\relax\aftergroup\noexpand
  \else\ifcat\noexpand#1\noexpand~\aftergroup\noexpand\fi\fi
  \aftergroup#1\sanitize}
\def\endsanity\endsanity{}

\def\note#1#2{\hbox to2in{\strut#1\quad\dotfill\quad#2}}
\def\boxit#1{\setbox4=\hbox{\kern1pt#1\kern1pt}
  \hbox{\vrule\vbox{\hrule\kern1pt\box4\kern1pt\hrule}\vrule}}
\def\halmos{\hbox{\am\char'3}} 
\def\qed#1\par{\message{.                                }\setbox1=\hbox{#1}%
  \ifdim\wd1>0pt\setbox\qedbox=\hbox{\halmos\raise2pt\hbox{\fiverm #1}}\fi
  \kern5pt\lower 2pt\hbox{\box\qedbox}\proofskip\goodbreak\endgroup}

\def\sectionsymbol{\S}
\def\k{\kappa}
\def\g{\gamma}
\def\a{\alpha}
\def\b{\beta}
\def\d{\delta}
\def\s{\sigma}
\def\t{\tau}
\def\l{\lambda}
\def\z{\zeta}
\def\I1{\mathop{\hbox{\sc i}_1}}
\def\w{\omega}
\def\P{{\mathchoice{\hbox{\bm P}}{\hbox{\bm P}}
         {\hbox{\tenbm P}}{\hbox{\sevenbm P}}}}
\def\Q{{\mathchoice{\hbox{\bm Q}}{\hbox{\bm Q}}
         {\hbox{\tenbm Q}}{\hbox{\sevenbm Q}}}}

\def\X{{\mathchoice{\hbox{\bm X}}{\hbox{\bm X}}
         {\hbox{\tenbm X}}{\hbox{\sevenbm X}}}}
\def\card#1{\left|#1\right|}

\def\id{\mathop{\hbox{\tenrm id}}}

\def\elesub{\prec}

\def\unifto{\buildrel\lower 7pt\hbox{$\to$}\over\to}

\def\cof{\mathop{\rm cof}\nolimits}
\def\cp{\mathop{\rm cp}\nolimits}

\def\ORD{\hbox{\sc ord}}

\def\CH{\hbox{\sc ch}}

\def\plus{^{\scriptscriptstyle +}}

\def\in{\mathrel{\mathchoice{\raise 
1pt\hbox{$\scriptstyle\cal\char'62$}}
         {\raise 1pt\hbox{$\scriptstyle\cal\char'62$}}
         {\raise .5pt\hbox{$\scriptscriptstyle\cal\char'62$}}
         {\hbox{$\scriptscriptstyle\cal\char'62$}}}\penalty700{}}
\def\ni{\mathrel{\mathchoice{\raise 1pt\hbox{$\scriptstyle\cal\char'63$}}
                   {\raise 1pt\hbox{$\scriptstyle\cal\char'63$}}
                   {\raise .5pt\hbox{$\scriptscriptstyle\cal\char'63$}}
                   {\hbox{$\scriptscriptstyle\cal\char'63$}}}\penalty700}
\def\of{\mathrel{\mathchoice{\raise 1pt\hbox{$\scriptstyle\subseteq$}}
                   {\raise 1pt\hbox{$\scriptstyle\subseteq$}}
                   {\raise .5pt\hbox{$\scriptscriptstyle\subseteq$}}
                   {\hbox{$\scriptscriptstyle\subseteq$}}}}
\def\fo{\mathrel{\mathchoice{\raise 1pt\hbox{$\scriptstyle\supseteq$}}
                   {\raise 1pt\hbox{$\scriptstyle\supseteq$}}
                   {\raise .5pt\hbox{$\scriptscriptstyle\supseteq$}}
                   {\hbox{$\scriptscriptstyle\supseteq$}}}}
\def\notin{\mathrel{\mathchoice
  {\raise 1pt\hbox{\rlap{$\scriptstyle\;|$}$\scriptstyle\cal\char'62$}}
  {\raise 1pt\hbox{\rlap{$\scriptstyle\kern2pt 
          |$}$\scriptstyle\cal\char'62$}}
  {\raise .5pt\hbox{\rlap{$\scriptscriptstyle\, |$}$\scriptscriptstyle
      \cal\char'62$}}
  {\hbox{\rlap{$\scriptscriptstyle\, |$}$\scriptscriptstyle
     \cal\char'62$}}}%
  \penalty700}

\def\ofnoteq{\mathbin{\hbox{\bm\char'050}}}

\def\and{\mathrel{\kern1pt\&\kern1pt}}
\def\iff{\mathrel{\leftrightarrow}}

\def\implies{\rightarrow}

\def\union{\cup}

\def\intersect{\cap}

\def\ot{\mathop{\rm ot}\nolimits}

\def\cross{\times}

\def\[#1]{\left[\vphantom{\bigm|}#1\right]}
\def\<#1>{\langle\,#1\,\rangle}

\def\sat{\models}

\def\image{\mathbin{\hbox{\tt\char'42}}}
\def\restrict{\mathbin{\mathchoice{\hbox{\am\char'26}}{\hbox{\am\char'26}}{\hbox{\eightam\char'26}}{\hbox{\sixam\char'26}}}}
\def\force{\mathbin{\hbox{\am\char'15}}}

\def\emptyset{\mathord{\hbox{\bm\char'77}}}

\def\concat{\mathbin{{\,\hat{ }\,}}}

\def\st{\mid}
\def\seq<#1>{{\def\st{\mid\penalty650}\left<\,#1\,\right>}}

\def\set#1{\{\,#1\,\}}

\def\th{{\hbox{\fiverm th}}}

\def\forces{\force}
\def\lttheta{{\raise 1pt\hbox{$\scriptstyle<$}\theta}}

\def\I1{\mathop{\hbox{\sc i}_1}}
\def\ltk{{{\scriptstyle<}\k}}
\def\ltl{{{\scriptstyle<}\l}}
\def\ltg{{{\scriptstyle<}\g}}

\def\lteb{{{\scriptstyle\leq}\b}}
\def\lted{{{\scriptstyle\leq}\d}}

\def\Qdot{\dot\Q}

\def\Adot{\dot A}
\def\Qdot{\dot\Q}
\def\qdot{\dot q}

\def\hdot{\dot h}
\def\rdot{\dot r}
\def\sdot{\dot s}

\def\jmu{j_\mu}
\def\jnu{j_\nu}

\def\ltg{{\scriptscriptstyle<}\g}

\def\qdot{\dot q}
\def\rdot{\dot r}
\def\sdot{\dot s}

\def\Vbar{{\overline V}}
\def\Mbar{{\overline M}}

\font\arrow=line10 scaled \magstep1
\def\makeline#1.{\hbox{\arrow\char#1}}
\def\makearrow#1.#2.{\hbox{\arrow\char#1\llap{\char#2}}}
\def\definelinesandarrows#1.#2.#3.#4.#5.{
   \expandafter\edef\csname#4line\endcsname{\makeline#1.}
   \expandafter\edef\csname#4arrow\endcsname{\makearrow#1.#2.}
   \expandafter\edef\csname#5line\endcsname{\makeline#1.}
   \expandafter\edef\csname#5arrow\endcsname{\makearrow#1.#3.}}
\definelinesandarrows 0.18.9.ne.sw.
\definelinesandarrows 1.21.11.nnne.sssw.
\definelinesandarrows 2.14.13.nnnne.ssssw.
\definelinesandarrows 3.23.15.nnnnne.sssssw.
\definelinesandarrows 4.23.15.nnnnnne.ssssssw.
\definelinesandarrows 10.30.29.nne.ssw.
\definelinesandarrows 16.49.41.neeeeee.swwwwww.
\definelinesandarrows 17.51.43.neeee.swwww.
\definelinesandarrows 19.55.47.nehuh.swhuh.
\definelinesandarrows 24.58.41.neeeeeee.swwwwwww.
\definelinesandarrows 26.62.9.neee.swww.
\definelinesandarrows 33.49.25.neeeee.swwwww.
\definelinesandarrows 35.62.61.nee.sww.
\definelinesandarrows 64.82.73.se.nw.
\definelinesandarrows 65.85.75.ssse.nnnw.
\definelinesandarrows 66.78.77.sssse.nnnnw.
\definelinesandarrows 67.87.79.ssssse.nnnnnw.
\definelinesandarrows 68.87.79.sssssse.nnnnnnw.
\definelinesandarrows 74.94.93.sse.nnw.
\definelinesandarrows 80.113.105.seeeeee.nwwwwww.
\definelinesandarrows 81.115.107.seeee.nwwww.
\definelinesandarrows 99.126.125.see.nww.
\def\sejoin#1#2{\setbox1=\hbox{#1}\setbox2=\hbox{#2}%
  \hbox{\vbox{\hbox{\copy1\kern\wd2}\nointerlineskip
              \hbox{\kern\wd1\box2}}}}
\def\nejoin#1#2{\setbox1=\hbox{#1}\setbox2=\hbox{#2}%
  \hbox{\vbox{\hbox{\kern\wd1\copy2}\nointerlineskip\hbox{\copy1\kern\wd2}}}}
\newdimen\hnudge
\newdimen\vnudge
\newdimen\hnudgedefault
\newdimen\vnudgedefault

\def\SEdefaultnudge{\hnudge=-16pt\vnudge=20pt}
\def\Edefaultnudge{\hnudge=-25pt\vnudge=6pt}

\def\longEdefaultnudge{\hnudge=-5pt\vnudge=6pt}
\def\nudgeright#1pt{\advance\hnudge by#1pt}
\def\nudgeleft#1pt{\advance\hnudge by-#1pt}
\def\nudgeup#1pt{\advance\vnudge by#1pt}
\def\nudgedown#1pt{\advance\vnudge by-#1pt}
\def\label#1{\smash{\llap{\kern\hnudge
                   \raise\vnudge\rlap{$\scriptstyle#1$}\hfill}}}

\def\SEarrow{\SEdefaultnudge
             \sejoin\seeline{\sejoin\seeline{\sejoin\seeline\seearrow}}}

\def\Earrow{\Edefaultnudge\setbox1=\hbox{\SEarrow}
 \hbox{\raise 2pt\hbox{\vrule height-.4pt depth.8ptwidth\wd1\kern2pt
       \llap{\arrow\char'55}}}}
\def\longEarrow{\longEdefaultnudge\setbox1=\hbox{\SEarrow}
      \rlap{\hskip-1.25\wd1\raise 2pt
            \hbox{\vrule height-.4pt depth.8ptwidth2.5\wd1\kern2pt
            \llap{\arrow\char'55}}}}

\looselineskip
\def\Mbar{{\overline M}}
\def\Vbar{{\overline V}}
\def\tdot{\dot t}
\centerline{[submitted to the Journal of Mathematical Logic, January, 1999]}

\title Gap Forcing

\author Joel David Hamkins\cr
            City University of New York\cr
            {\tentt http://www.math.csi.cuny.edu/$\sim$hamkins}\cr

\abstract. In this paper I generalize the landmark Levy-Solovay Theorem \cite[LevSol67], which limits the kind of large cardinal embeddings that can exist in a small forcing extension, to a broad new class of forcing notions, a class that includes many of the forcing iterations most commonly found in the large cardinal literature. The fact is that after such forcing, every embedding satisfying a mild closure requirement lifts an embedding from the ground model. A consequence is that such forcing can create no new weakly compact cardinals, measurable cardinals, strong cardinals, Woodin cardinals, strongly compact cardinals, supercompact cardinals, almost huge cardinals, or huge cardinals, and so on. 

\acknowledgement My research has been supported in part by grants from the PSC-CUNY Research Foundation and from the Japan Society for the Promotion of Science. I would like to thank my gracious hosts at Kobe University in Japan for their generous hospitality. This paper follows up an earlier announcement of the main theorem appearing, without technical details, in \cite[Ham$\infty$]. 

Small forcing in a large cardinal context, that is, forcing with a poset $\P$ of cardinality less than whatever large cardinal $\k$ is under consideration, is today generally looked upon as benign. This outlook is largely due to the landmark Levy-Solovay theorem \cite[LevSol67], which asserts that small forcing does not affect the measurability of any cardinal. (Specifically, the theorem says that if a forcing notion $\P$ has size less than $\k$, then the ground model $V$ and the forcing extension $V^\P$ agree on the measurability of $\k$ in a strong way: all the ground model measures on $\k$ generate as filters measures in the forcing extension, the corresponding ultrapower embeddings lift uniquely from the ground model to the forcing extension and all the measures and ultrapower embeddings in the forcing extension arise in this way.) Since the Levy-Solovay argument generalizes to the other large cardinals whose existence is witnessed by certain kinds of measures or ultrapowers, such as strongly compact cardinals, supercompact cardinals, almost huge cardinals and so on, one is led to the broad conclusion that small forcing is harmless; one can understand the measures in a small forcing extension by their relation to the measures existing already in the ground model. 

Historically, the Levy-Solovay theorem addressed G\"odel's hope that large cardinals would settle the Continuum Hypothesis (\CH). G\"odel, encouraged by Scott's \cite[Sco61] theorem showing that the existence of a measurable cardinal implies $V\not=L$, had hoped that large cardinals would settle the \CH\ in the negative. But since one can force the \CH\ to hold or fail quite easily with small forcing, the conclusion is inescapable that large cardinals simply have no bearing whatsoever on the Continuum Hypothesis.

Since that time, set theorists have developed sophisticated tools to combine the two central set theoretic topics of forcing and large cardinals. The usual procedure when forcing with a large cardinal $\k$ whose largeness is witnessed by the existence of a certain kind of elementary embedding $j:V\to M$ is to lift the embedding to the forcing extension $j:V[G]\to M[j(G)]$ and argue that this lifted embedding witnesses that $\k$ retains the desired large cardinal property in $V[G]$. In this way, one is led to consider how a measure $\mu$ in the ground model $V$ can relate to a measure $\nu$ in the forcing extension $V[G]$. The measure $\mu$ may {\df extend} to $\nu$ in the simple sense that $\mu\of\nu$ or it may {\df lift} to $\nu$ when the ground model ultrapower $\jmu:V\to M$ agrees with the larger ultrapower $\jnu:V[G]\to M[\jnu(G)]$ on the common domain $V$. Expressed in this terminology, the Levy-Solovay theorem asserts that after small forcing every measure in the ground model both lifts and extends to a measure in the forcing extension and, conversely, every measure in the extension both lifts and extends a measure in the ground model. 

The truth, however, is that in the large cardinal context most small forcing is, as it were, too small. Rather, one often wants to perform long iterations going up to and often beyond the large cardinal $\k$ in question. With a supercompact cardinal $\k$, for example, one often sees reverse Easton $\k$-iterations along the lines of Silver forcing \cite[Sil71] or the Laver preparation \cite[Lav78]. What we would really like is a generalization of the Levy-Solovay theorem that would allow us to understand and control the sorts of embeddings and measures added by these more powerful and useful forcing notions.

Here, I prove such a generalization. For a vast class of forcing notions, including the iterations I have just mentioned, the fact is that every embedding $j:V[G]\to M[j(G)]$ in the extension that satisfies a mild closure condition lifts an embedding $j:V\to M$ from the ground model. In particular, every measure in $V[G]$ concentrating on a set in $V$ extends a measure on that set in $V$. From this general fact, I deduce that forcing of this type creates no new weakly compact cardinals, measurable cardinals, strong cardinals, Woodin cardinals, supercompact cardinals, or huge cardinals and so on. 

The class of forcing notions for which the theorem applies is quite broad. All that is required is that the forcing admit a {\df gap} at some $\d$ below the cardinal $\k$ in question in the sense that the forcing factors as $\P*\Qdot$ where $\P$ is nontrivial, $\card{\P}<\d$ and $\forces\Qdot$ is $\lted$-strategically closed. (A forcing notion is $\lted$-strategically closed when the second player has a strategy enabling her to survive through all the limits in the game in which the players alternately play conditions to build a descending $(\d+1)$-sequence through the poset, with the second player playing at limit stages.) The Laver preparation, for example, admits a gap between any two stages of forcing. Indeed, in the Laver preparation, the tail forcing is fully directed closed, not merely closed or strategically closed. And the same holds for many of the other reverse Easton iterations one commonly finds in the literature. Moreover, in practice one can often simply preface whatever strategically closed forcing is at hand with some harmless small forcing, such as the forcing to add a single Cohen real, and thereby introduce a gap at $\d=\omega_1$. Further, because $\Qdot$ can be trivial, gap forcing includes all small forcing notions. Examples of useful gap forcing notions are abundant. 

An embedding $j:\Vbar\to \Mbar$ is {\df amenable} to $\Vbar$ when $j\restrict A\in\Vbar$ for any $A\in\Vbar$. 

\quiet\theorem Gap Forcing Theorem. Suppose that $V[G]$ is a forcing extension obtained by forcing that admits a gap at some $\d$ below $\k$ and $j:V[G]\to M[j(G)]$ is an embedding with critical point $\k$ for which $M[j(G)]\of V[G]$ and $M[j(G)]^\d\of M[j(G)]$ in $V[G]$. Then $M\of V$; indeed $M=V\intersect M[j(G)]$. If the full embedding $j$ is amenable to $V[G]$, then the restricted embedding $j\restrict V:V\to M$ is amenable to $V$. And if $j$ is definable from parameters (such as a measure or extender) in $V[G]$, then the restricted embedding $j\restrict V$ is definable from the names of those parameters in $V$.

A weaker precursor to this theorem appeared in \cite[Ham98a]\footnote{${}^1$}{The current proof addresses what is probably an inadequate discussion of $\ddot s$ in that proof.}. The Gap Forcing Theorem here answers all of the open questions asked in \cite[Ham98a] and establishes a strong generalization of the Gap Forcing Conjecture of that paper, which asserted that after forcing with a very low gap every supercompactness embedding is the lift of an embedding from the ground model. The current theorem implies much more: any kind of ultrapower embedding is a lift. 

In order to avoid confusion on a subtle point, let me remark that given any embedding $j:V[G]\to\Mbar$ we can let $M=\union\set{j(V_\a)\st\a\in\ORD}$, and it is not difficult to see that $j(G)$ is $M$-generic, that $\Mbar=M[j(G)]$ and moreover that $j\restrict V:V\to M$. Thus, while the statement of the theorem concerns embeddings of the form $j:V[G]\to M[j(G)]$, this form of embedding is fully general. 

For those readers who are not completely familiar with the bizarre sorts of embeddings $j:V[G]\to M[j(G)]$ that can exist in a forcing extension, let me stress that in general, quite apart from the question of whether $j$ lifts an embedding from the ground model, one must not presume even that $M\of V$. For example, if $\k$ is a Laver indestructible supercompact cardinal in $V$ and we force to add a Cohen subset $A\of\k$ (by itself, this forcing does not admit a gap below $\k$), then $\k$ remains supercompact in the extension $V[A]$, but any embedding $j:V[A]\to M[j(A)]$ must have $A\in M$ and therefore $M\not\of V$. The point is that the theorem really does identify a serious, useful limitation on the sorts of embeddings that exist in a gap forcing extension. 

My proof will proceed through a sequence of lemmas. A variation of the Key Lemma first appeared in \cite[Ham98a] and \cite[Ham98b] and was subsequently modified and appealed to in \cite[HamShl98], but to be thorough I include the proof here. Other important techniques are adapted from Woodin's proof of the Levy-Solovay Theorem for strong cardinals (see \cite[HamWdn]); indeed, Woodin's techniques are peppered amongst the proofs of several of the lemmas below, and I could not have proved the theorem without them.

Let me define that a sequence in a forcing extension is {\df fresh} when it is not in the ground model but all of its proper initial segments are. Thus, it is a new path through a tree in the ground model.

\lemma Key Lemma. If $\card{\P}\leq\b$, $\forces\Qdot$ is $\lteb$-strategically closed and $\cof(\theta)>\b$, then $\P*\Qdot$ adds no fresh $\theta$-sequences. 

\proof It suffices to consider only sequences of ordinals. Furthermore, since any fresh $\theta$-sequence of ordinals below $\xi$ may be easily coded with a fresh binary sequence of ordinal length $\xi\cdot\theta$, which has the same cofinality as $\theta$, it suffices to prove only that no fresh binary sequences are added. So, suppose towards a contradiction that $\t$ is the $\P*\Qdot$-name of a fresh binary $\theta$-sequence, so that 
$$\forces_{\P*\Qdot}\,\t\in 2^{\check\theta}\and\t\notin\check V\and\forall\l<\check\theta\,(\t\restrict\l\in\check V).$$
Since $\P$ is nontrivial, by refinining below a condition if necessary we may assume it adds a new subset of some minimal $\g\leq\b$, so that for some name $\hdot$:
$$\forces_{\P}\,\hdot\in 2^{\check\g}\and\hdot\notin\check V\and\forall\a<\check\g\,(\hdot\restrict\a\in\check V).$$
For every condition $\<p,\qdot>\in\P*\Qdot$ let $b_{\<p,\qdot>}$ be the longest sequence $b$ such that $\<p,\qdot>\forces\check b\of\t$. Note that $\cof(\theta)>\b$ is preserved by both $\P$ and $\Q$. 

I claim that a certain weak Prikry property holds, namely, that there is a condition $\<p,\qdot>$ such that for for any $\l<\theta$ and any stronger condition of the form $\<p,\rdot>$ there is an even stronger condition of the form $\<p,\sdot>$ that decides $\t\restrict\l$. That is, below $\<p,\qdot>$ the first coordinate need not change in order to decide more and more of $\t$. To see why this is so, suppose $g*G$ is $V$-generic for $\P*\Qdot$. For every $\l<\theta$ there is a condition $\<p_\l,\qdot_\l>\in g*G$ that decides $\t\restrict\l$. Since $\cof(\theta)>\b$, it must be that a single condition $p$ is used for unboundedly many $p_\l$. Thus, in fact, this condition $p$ could have been used for every $\l$. So for every $\l$ there is a name $\qdot$ such that $\<p,\qdot>\in g*G$ decides $\t\restrict\l$. By strengthening $p$ if necessary, we may suppose that this state of affairs is forced by a condition of the form $\<p,\qdot>$. What this means is that for any $\l$ and any stronger $\<p,\rdot>$ there is an even stronger $\<p,\sdot>$ that decides $\t\restrict\l$, as I claimed.

Since no condition decides all of $\t$, it follows from this that for any condition $\<p,\rdot>\leq\<p,\qdot>$ there are names $\rdot_0$ and $\rdot_1$ such that $\<p,\rdot_0>$, $\<p,\rdot_1>\leq\<p,\qdot>$ and $b_{\<p,\rdot_0>}\perp b_{\<p,\rdot_1>}$. 

Now I will iterate this fact by constructing in $V$ a binary branching tree whose paths represent (names for) the first player's plays in the game corresponding to $\Qdot$. Using a name $\dot\s$ that with full boolean value names a strategy witnessing that $\Qdot$ is $\lteb$-strategically closed in $V^\P$, the basic picture is that while the second player obeys $\dot\s$, the tree will branch for the first player with moves corresponding to the conditions $\rdot_i$ given by the previous paragraph. Specifically, I will assign in $V$ to each $t\in 2^{{<}\g}$ a name $\qdot_t$ so that along any branch in $2^\g$ the condition $p$ forces that the names give rise to the first player's moves in a play through $\Qdot$ that accords with the strategy $\dot\s$. That is, the next move is always below $\dot\s$ of the previous moves. The first player begins with $\qdot_{\emptyset}=\qdot$. If $\qdot_t$ is defined, let $\rdot_t$ be the name of the condition obtained by applying the strategy against the play up to this point, i.e. the play in which the first player plays $\qdot_s$ for $s\of t$. By induction, $p$ forces that these conditions give rise to a play according to $\dot\s$, and so $p$ forces that $\rdot_t$ is stronger than all $\qdot_s$ for $s\of t$. Now, by the previous paragraph the first player may reply with either $\qdot_{t\concat 0}$ or $\qdot_{t\concat 1}$ chosen so that $\<p,\qdot_{t\concat i}>\leq\<p,\rdot_t>$ and $b_{\<p,\qdot_{t\concat 0}>}\perp b_{\<p,\qdot_{t\concat 1}>}$. Similarly, if $t$ has limit ordinal length, then since the strategy is forced to be winning for the second player, there will be a condition $\rdot_t$ that is the result of the strategy $\dot\s$ applied to the previous play $\<\qdot_s\st s\ofnoteq t>$, and we may therefore have the first player choose any $\qdot_t$ such that $\<p,\qdot_t>\leq\<p,\rdot_t>$ in order to continue the iteration. The effects of this construction are first, that whenever $t\of\bar t$, then $\<p,\qdot_{\bar t}>\leq\<p,\qdot_t>$, and second, that $b_{\<p,\qdot_{t\concat 0}>}\perp b_{\<p,\qdot_{t\concat 1}>}$. The map $t\mapsto\qdot_t$ lies in $V$. 

Now suppose that $g*G$ is $V$-generic below the condition $\<p,\qdot>$. In $V[g]$, let $h=\hdot_g$ be the new $\g$-sequence added by $\P$; let $q_t=(\qdot_t)_g$ be the interpretation of the names constructed in the previous paragraph; and let $\s=(\dot\s)_g$ be the interpretation of the strategy. By the assumption on $\hdot$, every initial segment $t\ofnoteq h$ lies in $V$. By construction, the sequence $\<q_t\st t\ofnoteq h>$ represents the plays of the first player in a play that accords with the strategy $\s$. Thus, since the strategy is winning for the second player, there is a condition $r$ below all of them (i.e. the $\g^\th$ move). Thus, $r$ forces that $b=\union_{t\ofnoteq h}b_{\<p,\qdot_t>}$ is a proper initial segment of $\t$, and consequently $b\in V$. By construction, however, for any $t\in 2^{\ltg}$ we know $t\of h$ exactly when $b_{\<p,\qdot_t>}\of b$, since whenever $t\concat i$ first deviates from $h$ the construction ensures that $b_{\<p,\qdot_{t\concat i}>}$ deviates from $b$. We conclude that $h\in V$, a contradiction.\qed 

Let me now continue with the proof of the theorem. Suppose that $V[G]$ is a forcing extension obtained by forcing that admits a gap at $\d<\k$ and $j:V[G]\to M[j(G)]$ is an embedding with critical point $\k$ such that $M[j(G)]\of V[G]$ and $M[j(G)]^\d\of M[j(G)]$ in $V[G]$. Exhibiting the gap, we have $V[G]=V[g][H]$ where $g*H\of\P*\Qdot$ is $V$-generic for nontrivial forcing $\P$ with $\card{\P}<\d$ and $\forces\Qdot$ is $\lted$-strategically closed. The embedding can therefore be written as $j:V[g][H]\to M[g][j(H)]$. I may assume that $\d$ is regular, since it might as well be $\card{\P}\plus$. Since the critical point of $j$ is $\k$, every set in $V_\k$ is fixed by $j$. It follows that $V_\k=M_\k$. In the next few lemmas, I will show even more agreement between $M$ and $V$. 

\lemma. Every set of ordinals $\s$ in $V[G]$ of size $\d$ is covered by a set $\t$ in $M\intersect V$ of size $\d$. 

\proof Since $\s$ has size $\d$, it must be in both $V[g]$ and $M[g]$. Thus, using the names in $V$ and $M$, there are sets $s_0\in V$ and $s_1\in M$ of size $\d$ such that $\s\of s_0$ and $\s\of s_1$. Iterating this idea, bouncing between sets in $M$ and sets in $V$, we can build in $V[G]$ an increasing sequence of sets $\vec\s=\<\s_\a\st \a<\d>$ such that $\s_0=\s$, $\a<\b\implies\s_\a\of\s_\b$, and for cofinally many $\a$, $\s_\a\in V$ and for cofinally many $\a$, $\s_\a\in M$. Let $\t=\union\vec\s$. Thus certainly $\s\of\t$ and $\t$ has size $\d$. It remains to show $\t\in M\intersect V$. By the strategic closure of $\Q$ we know $\vec\s\in V[g]$, and so it has a name $\sdot\in V$. Since cofinally often $\s_\a\in V$, there must be conditions in $g$ forcing each instance of this, but since $\card{\P}<\d$ and $\d$ is regular, a single condition $p\in g$ must work unboundedly often, and decide unboundedly many elements of $\sdot$. Thus, $p$ also decides the union, and so $\t\in V$. Similarly, by the closure of the embedding it must be that $\vec\s\in M[j(G)]$ and consequently by the strategic closure of $j(\Q)$ actually $\vec\s\in M[g]$. Thus, it has a name $\tdot\in M$, and again because cofinally often $\s_\a\in M$ there must be a single condition $p\in g$ deciding unboundedly many many elements of $\tdot$. Thus, this condition decides the union, and so $\t\in M$, as desired.\qed

\lemma. $M$ and $V$ have the same $\d$-sequences of ordinals. 

\proof It suffices to show that $[\ORD]^\d$ is the same in $M$ and $V$. Suppose that $\s\of\ORD$ has size $\d$ and $\s$ is in either $M$ or $V$. By the previous lemma there is a set $\t\in V\intersect M$ of size $\d$ such that $\s\of\t$. In both $M$ and $V$ we may enumerate $\t=\set{\b_\a\st \a<\g}$ in increasing order, where $\g=\ot(\t)<\d\plus$. Let $A=\set{\a\st \b_\a\in\s}$. This set is definable from $\s$ and $\t$ and therefore must be in either $M$ or $V$, respectively, as $\s$ is in either $M$ or $V$. But since $A\of\g$, it must be in $V_\k=M_\k$, and so it is in {\it both} $M$ and $V$. Thus, $\s=\set{\b_\a\st \s\in A}$ is also in both $M$ and $V$, as desired.\qed

\lemma. $M\of V$. 

\proof It suffices to show that $P(\theta)^M\of V$ for every ordinal $\theta$. Suppose $A\of\theta$ is in $M$. By induction, I may assume that every initial segment of $A$ is in $V$. If $\cof(\theta)\geq\d$ then $A$ must itself be in $V$, for otherwise it would be fresh over $V$, contradicting the Key Lemma. So we may assume that $\cof(\theta)<\d$. Thus, by the distributivity of $\Q$, it follows that $A\in V[g]$, and so $A=\Adot_g$ for some name $\Adot\in V$. Pick some enormous $\z$ and an elementary substructure $\X\elesub V_\z$ of size $\d$ containing $\Adot$ and $\P$ as well as every element of $\P$. It follows that $g$ is $\X$-generic, that $\X[g]\elesub V_\z[g]$ and furthermore that $\X$ and $\X[g]$ have the same ordinals. Since $\X\intersect\ORD$ is a set of ordinals in $V$ of size $\d$, by the previous lemma it must also be in $M$. And since $A\in M$, it follows that $a=A\intersect \X$ is also in $M$, and so again by the previous lemma, $a$ is in $V$. Thus, there is some condition $p\in g$ that forces $\X\intersect\Adot=\check a$. That is to say, $p$ decides $\Adot(\a)$ for every $\a\in\X$. Thus, 
$$\X\sat\forall\a\hbox{($p$ decides $\Adot(\a)$)}.$$
By elementarity, it must be that $V_\z$ also satisfies this, and so $p$ decides $\Adot(\a)$ for all $\a$. Thus, $A\in V$, as desired.\qed

By simply enumerating any set in $M$, it follows from the previous two lemmas that $M^\d\of M$ in $V$. 

\lemma. $M=V\intersect M[j(G)]$. 

\proof Since $M\of V$ by the previous lemma it follows that $M\of V\intersect M[j(G)]$. For the converse, let me first show that any set of ordinals $A$ in both $V$ and $M[j(G)]$ is in $M$. Suppose $A\of\theta$ and, by induction, all the proper initial segments of $A$ are in $M$. If the cofinality of $\theta$ is at least $\d$, then $A$ must be in $M$ because otherwise it would be a fresh set added by $j(G)$ over $M$ in violation of the Key Lemma. So I may assume that $\cof(\theta)<\d$. Thus, $A$ is the union of a $\d$-sequence of sets in $M$. Since $M^\d\of M$ in $V$, it follows that $A$ is in $M$, as desired. Suppose now that $a$ is an arbitrary set in $V\intersect M[j(G)]$ and, by $\in$-induction, that every element of $a$ is in $M$. It follows that $a$ is a subset of an element $b\in M$. Enumerate $b=\set{b_\a\st \a<\theta}$ in $M$ and observe that it is enough to know that the set $A=\set{\a<\theta\st b_\a\in a}$ is in $M$. But this is a set ordinals in $V\intersect M[j(G)]$, and so the proof is complete.\qed

\lemma. If the full embedding $j:V[G]\to M[j(G)]$ is amenable to $V[G]$, then the restricted embedding $j\restrict V:V\to M$ is amenable to $V$. 

\proof Suppose that $j:V[G]\to M[j(G)]$ is amenable to $V[G]$. In order to show that $j\restrict V:V\to M$ is amenable to $V$, I must show that $j\restrict A\in V$ for any $A\in V$. Using enumerations of the sets in $V$, it suffices to show that $j\restrict\theta\in V$ for every ordinal $\theta$. And to prove this, it suffices to show that $j\image\theta\in V$ for every ordinal $\theta$. Let $A=j\image\theta$, and suppose by induction that every initial segment of $A$ is in $V$. By the amenability of the full embedding, we know that $A\in V[G]$. If $\cof(\theta)\geq\d$ then $A$ must be in $V$ for otherwise it would be fresh over $V$, in violation of the Key Lemma. So I may assume that $\cof(\theta)<\d$. Consequently, by the distributivity of $\Q$, it must be that $A\in V[g]$, and so $A=\Adot_g$ for some name $\Adot\in V$. Again choose some large $\z$ and $\X\elesub V_\z$ of size $\d$ containing $\Adot$ and $\P$ as well as every element of $\P$. It follows that $\X\intersect\ORD=\X[g]\intersect\ORD$. The set $\X\intersect\ORD$ is a set of ordinals of size $\d$ in $V$, and consequently it is in $M$ by the lemma above. Let $a=A\intersect\X=A\intersect\X[g]$. Since this is a subset of $j\image\theta$ of size $\d<\k$, it must be equal to $j\image b=j(b)$ for some set $b\of\theta$ of size $\d$. By the cover lemma above, there is a set $c$ in both $M$ and $V$ such that $b\of c$ and $c$ has size $\d$. Now simply compute $a=j\image b\of (j\image c)\intersect\X\of (j\image\theta)\intersect\X=a$, and so $a=(j\image c)\intersect\X$. But $j\image c=j(c)\in M\of V$, and so $a$ is in $V$. Now, continuing as in the previous lemma, there must be a condition $p\in g$ forcing this. So $p$ decides $\Adot(\a)$ for every $\a\in\X$. By the elementarity of $\X\elesub V_\z$ it must be that $p$ decides $\Adot(\a)$ for every ordinal $\a$. Thus, $A$ is in $V$, as desired.\qed

\lemma. If the full embedding $j:V[G]\to M[j(G)]$ is definable from parameters (such as a measure or extender) in $V[G]$, then the restricted embedding $j\restrict V:V\to M$ is definable from the names of those parameters in $V$.

\proof This follows actually from the previous lemma. Suppose that $j:V[G]\to M[j(G)]$ is definable from the parameter $z$ in $V[G]$. Thus, there is some formula $\varphi$  such that $j(a)=b$ exactly when $V[G]\sat\varphi[a,b,z]$. Fix a name $\dot z$ for $z$. Thus, for $a\in V$ we have $j(a)=b$ exactly when some $p\in G$ forces that $\varphi(\check a,\check b,\dot z)$. Since any definable embedding is amenable, it follows by the previous lemma that $j\restrict V_\theta$ is in $V$ for every $\theta$. Thus, for every $\theta$ there is a condition $p\in G$ that forces that the relation $\varphi(\check a,\check b,\dot z)$ in $V[G]$ agrees with the relation $j(a)=b$ for $a\in V_\theta$ and $b\in V_{j(\theta)}$. That is, $p$ forces that the relation $\varphi(\check a,\check b,\dot z)$ for $a$ and $b$ in the appropriate domain produces exactly the set $j\restrict V_\theta$. By the Axiom of Replacement, there must be a single $p$ that works for unboundedly many $\theta$. Thus, for this $p$ we know that $j(a)=b$ exactly when $p$ forces $\varphi(\check a,\check b,\dot z)$, for $a$ and $b$ in $V$. So $j\restrict V$ is definable from $\dot z$ in $V$.\qed

This completes the proof of the theorem. I will nevertheless quickly prove one additional lemma that will assist in the proofs of the corollaries to come. 

\lemma. Under the hypothesis of the theorem, 
\points 1. If $M[j(G)]^\l\of M[j(G)]$ in $V[G]$ then $M^\l\of M$ in $V$.\cr
           2. If $V_\l\of M[j(G)]$ then $V_\l\of M$.\cr

\ref\Assist

\proof For 1, if $M[j(G)]^\l\of M[j(G)]$ in $V[G]$ then any $\l$-sequence of elements of $M$ that lies in $V$ must lie in $V\intersect M[j(G)]$, and consequently in $M$. For 2, if $V_\l\of M[j(G)]$ then $V_\l\of V\intersect M[j(G)]=M$.\qed

One must take care with strongness embeddings in order to satisfy the closure hypothesis in the theorem. A cardinal $\k$ is {\df $\l$-strong} when there is an embedding $j:V\to M$ with critical point $\k$ such that $V_\l\of M$ and $j(\k)>\l$. Let me define that an embedding $j:V\to M$ is {\df $\b$-closed} when $M^\b\of M$. The problem with strongness embeddings, of course, is that they need not satisfy any degree of closure. By factoring through by the canonical extender, however, one obtains a {\df natural} embedding, meaning in addition that $M=\set{j(h)(s)\st h\in V\and s\in V_\l}$. And for almost every $\l$ these natural embeddings do satisfy the closure hypothesis of the theorem.

\lemma. If $\k$ is $\l$-strong, then the natural $\l$-strongness embeddings $j:V \to M$ are $\k$-closed if $\l$ is a successor ordinal or a limit ordinal of cofinality above $\k$, and otherwise they are ${<}\cof(\l)$-closed. 
\ref\Strongness

\proof It suffices to consider $\l>\k$. Suppose that $j:V\to M$ is a natural $\l$-strongness embedding, so that $\cp(j)=\k$, $V_\l\of M$ and $M=\set{j(h)(s)\st s\in V_\l\and h\in V}$. In the first case, suppose that $\l=\xi+1$ and $\<j(h_\a)(s_\a)\st\a<\k>$ is a $\k$-sequence of elements from $M$. Since a $\k$-sequence of subsets of $V_\xi$ can be coded with a single subset of $V_\xi$, it follows that $\<s_\a\st\a<\k>$ is in $M$. Also, the sequence $\<j(h_\a)\st\a<\k>=j(\<h_\a\st\a<\k>)\restrict\k$ is in $M$. Thus, the sequence $\<j(h_\a)(s_\a)\st\a<\k>$ is in $M$, as desired. For the next case, when $\l$ is a limit ordinal of cofinality larger than $\k$, then on cofinality grounds the sequence $\<s_\a\st\a<\k>$ is in $V_\l$, and hence in $M$, so again $\<j(h_\a)(s_\a)\st\a<\k>$ is in $M$, as desired. Finally, suppose $\l$ is a limit ordinal and $\b<\cof(\l)\leq\k$. If $\<j(h_\a)(s_\a)\st\a<\b>$ is a sequence of elements from $M$, then again on cofinality grounds we know $\<s_\a\st\a<\b>$ is in $V_\l$ and hence in $M$, and so $\<j(h_\a)(s_\a)\st\a<\b>$ is in $M$, as desired.\qed

The consequence of this argument is that except for the limit ordinals of small cofinality, the Gap Forcing Theorem applies to strongness embeddings. 

I would like now to prove a series of corollaries to the Gap Forcing Theorem. I hope these corollaries tend to show that for a variety of large cardinals the restrictions identified in the theorem are severe. 

\corollary. Gap forcing creates no new weakly compact cardinals. If $\k$ is weakly compact after forcing with a gap below $\k$, then it was weakly compact in the ground model.

\proof Suppose that $\k$ is weakly compact in $V[G]$, a forcing extension obtained by forcing with a gap below $\k$. It follows that $\k$ is inaccessible in $V[G]$ and hence also in $V$. Thus, it remains only to prove that $\k$ has the tree property in $V$. If $T$ is a $\k$-tree in $V$, then by weak compactness it must have a $\k$-branch in $V[G]$. Since every initial segment of this branch is in $V$, it follows by the Key Lemma that the branch itself is in $V$, as desired.\qed

\corollary. After forcing with a gap below $\k$, every ultrapower embedding with critical point $\k$ in the extension lifts an embedding from the ground model, and every $\k$-complete measure in the extension that concentrates on a set in the ground model extends a measure in the ground model. 

\proof I am referring here not just to measures on $\k$, but to measures on an arbitrary set $D$, so that the corollary also covers the cases of, for example, supercompactness and hugeness measures. It is a standard fact that any ultrapower embedding $j:V[G]\to M[j(G)]$ by a measure $\mu$ in $V[G]$ is closed under $\k$-sequences where $\k=\cp(j)$. Since the forcing admits a gap below $\k$, the Gap Forcing Theorem implies that $j:V\to M$ is definable from parameters in $V$. If $\mu$ concentrates on some set $D\in V$, then since $X\in\mu\iff[\id]_\mu\in j(X)$, it follows that $\mu\intersect V\in V$ is a measure on $D$ in $V$, and the corollary is proved.\qed

\corollary. Gap forcing creates no new measurable cardinals. If $\k$ is measurable after forcing with a gap below $\k$, then $\k$ was measurable in the ground model and every measure on $\k$ in the extension extends a measure in the ground model. 

\proof This is a special case of the previous corollary.\qed

As a caution to the reader, let me stress that the corollary does not say that every ultrapower embedding $j:V[G]\to M[j(G)]$ in the extension is the lift of an ultrapower embedding $j\restrict V:V\to M$ in the ground model. Rather, one only knows that the restricted embedding $j\restrict V$ is definable from parameters in $V$. Indeed, it is possible to construct an example of a gap forcing extension $V[G]$ with an  embedding $j:V[G]\to M[j(G)]$ that is the ultrapower by a normal measure in $V[G]$ but the restriction $j\restrict V$ is not an ultrapower embedding at all, being instead some kind of strongness extender embedding. 

\corollary. Gap forcing creates no new strong cardinals. If $\k$ is $\l$-strong after forcing with a gap at $\d<\k$, and $\l$ is either a successor ordinal or has cofinality larger than $\d$, then $\k$ was $\l$-strong in the ground model. 

\proof Suppose that $V[G]$ is the forcing extension obtained by forcing with a gap below $\d$. Lemma \Strongness\ shows that if $\k$ is $\l$-strong in $V[G]$ for such a $\l$ as in the statement of the corollary, then there is an embedding $j:V[G]\to M[j(G)$ witnessing this that is closed under $\d$-sequences. Consequently, by the Gap Forcing Theorem, the restriction $j:V\to M$ is definable from parameters in $V$, and since $V_\l\of M[j(G)]$, Lemma \Assist\ implies that $V_\l\of M$. So $\k$ is $\l$-strong in $V$, as desired.\qed

What we actually have is the following:

\corollary. After forcing $\P$ of size less than $\d$, no further $\lted$-strategically closed forcing $\Q$ can increase the degree of strongness of any cardinal $\k>\d$. 

\proof Suppose that $\k$ is $\l$-strong in $V[g][H]$, the extension by $\P*\Qdot$, and $\k>\d$. In the first case, when $\l$ is either a successor ordinal or a limit ordinal of cofinality above $\d$, the previous corollary shows that $\k$ is $\l$-strong in $V$ and hence also in the small forcing extension $V[g]$. For the second, more difficult case, suppose that $\k$ is $\l$-strong in $V[g][H]$ and $\l$ is a limit ordinal with $\cof(\l)\leq\d$. Let $j:V[g][H]\to M[g][j(H)]$ be a $\l$-strong embedding by a canonical extender, so that $M[g][j(H)]=\set{j(h)(s)\st s\in V[g][H]_\l\and h\in V[g][H]}$. Thus, $j$ is the embedding induced by the extender $$E=\set{\<A,s>\st A\of V_\k\and s\in j(A)\and s\in V[g][H]_\l},$$ which is a subset of $P(V_\k)\cross V[g][H]_\l$. This extender is the union of the smaller extenders $E\restrict\b=E\intersect(P(V_\k)\cross V[g][H]_\b)$ for unboundedly many $\b<\l$. By the result of the previous corollary, we may assume that these smaller extenders each extend a strongness extender in $V$. Since each of these extenders extends uniquely to $V[g]$, the small forcing extension, it follows by the strategic closure of $\Qdot$ that $E\intersect V[g]$ is in $V[g]$ and hence $\k$ is $\l$-strong in $V[g]$, as desired.\qed

The two previous results are complicated somewhat by the intriguing possibility that small forcing could actually increase the degree of strongness of some cardinal. This question, an unresolved instance of the Levy-Solovay theorem, is raised in \cite[HamWdn]. One could ask the corresponding question replacing small forcing with gap forcing, {\it is it possible that forcing with a gap below $\k$ can increase the degree of strongness of $\k$?}\/ But the truth of the matter is that the previous corollary shows that if gap forcing $\P*\Qdot$ can increase the degree of strongness of a cardinal, then this increase is entirely due to the initial small forcing factor $\P$. And the only way this can occur is if a $\ltl$-strong cardinal is made $\l$-strong for some limit ordinal $\l$ of small cofinality.

\corollary. Gap forcing creates no new Woodin cardinals. If $\k$ is Woodin after forcing with a gap below $\k$, then $\k$ was Woodin in the ground model.

\proof If $\k$ is Woodin in $V[G]$, then for every $A\of\k$ there is a cardinal $\g<\k$ that is $\ltk$-strong for $A$, meaning that for every $\l<\k$ there is an embedding $j:V[G]\to M[j(G)]$ with critical point $\g$ such that $A\intersect\l=j(A)\intersect\l$. Such an embedding can be found that is $(\l+1)$-strong and induced by the canonical extender, so by Lemma \Strongness\ we may assume that $M[j(G)]$ is closed under $\g$-sequences. Further, such $\g$ must be unbounded in $\k$, so we may consider some such $\g$ above the gap in the forcing. Thus, for $A$ in the ground model, the Gap Forcing Theorem shows that the restricted embedding $j:V\to M$ witnesses the $\l$-strongness of $\g$ for $A$ in $V$, and so $\k$ was a Woodin cardinal in $V$, as desired.\qed 

Define that a forcing notion is {\df mild} relative to $\k$ when every set of ordinals of size less than $\k$ in the extension has a name of size less than $\k$ in the ground model. For example, the reverse Easton iterations one often finds in the literature are generally mild because the tail forcing is usually sufficiently distributive, and so any set of ordinals of size less than $\k$ is added by some stage before $\k$. Additionally, any $\k$-c.c. forcing is easily seen to be mild. 

\corollary. Mild gap forcing creates no new strongly compact cardinals. If $\k$ is $\l$-strongly compact after forcing that is mild relative to $\k$ and admits a gap below $\k$, then it was $\l$-strongly compact in the ground model; and every strong compactness measure in the extension is isomorphic to one that extends a strong compactness measure from the ground model.

\proof The point is that after mild forcing, every strong compactness measure $\mu$ on $P_\k\theta$ in the extension is isomorphic to a strong compactness measure $\tilde\mu$ that concentrates on $(P_\k\theta)^V$. To see why this is so, let $j:V[G]\to M[j(G)]$ be the ultrapower by $\mu$, and let $s=[\id]_\mu$. Thus, $j\image\theta\of s\of j(\theta)$ and $\card{s}<j(\k)$. By mildness $s$ has a name in $M$ of size less than $j(\k)$, and using this name we can construct a set $\tilde s\in M$ such that $j\image\theta\of\tilde s\of j(\theta)$ and $\card{\tilde s}<j(\k)$ in $M$. Furthermore, since $\mu$ is isomorphic to a measure concentrating on $\theta$, there must be some ordinal $\z<j(\theta)$ such that $M[j(G)]=\set{j(h)(\z)\st h\in V[G]}$. I may assume that the largest element of $\tilde s$ has the form $\<\a,\z>$, using a suitable definable pairing function, by simply adding such a point if necessary. Let $\tilde\mu$ be the measure germinated by $\tilde s$ via $j$, so that $X\in\tilde\mu\iff\tilde s\in j(X)$. Since $\tilde s$ is a subset of $j(\theta)$ of size less than $j(\k)$ in $M$, it follows that $\tilde\mu$ is a fine measure on $P_\k\theta$ in $V[G]$ that concentrates on $(P_\k\theta)^V$. I will now show that $\mu$ and $\tilde\mu$ are isomorphic. For this, it suffices by the seed theory of \cite[Ham97] to show that every element of $M[j(G)]$ is in the seed hull $\X=\set{j(h)(\tilde s)\st h\in V[G]}\elesub M[j(G)]$ of $\tilde s$. By the choice of $\tilde s$ we know that $\z\in\X$ and so it is easy to conclude that $j(h)(\z)\in\X$ for any function $h\in V[G]$, as desired. So every strong compactness measure is isomorphic to a strong compactness measure that concentrates on $(P_\k\theta)^V$. 

Now the corollary follows because the restricted embedding $j\restrict V:V\to M$ must be definable (with a name for $\mu$ as a parameter) in $V$ by the Gap Forcing Theorem, and using this embedding one can recover $\tilde\mu\intersect V$, which is easily seen to be a fine measure on $P_\k\theta$ in $V$, as desired.\qed

\corollary. Gap forcing creates no new supercompact cardinals. Indeed, it does not increase the degree of supercompactness of any cardinal. If $\k$ is $\l$-supercompact after forcing with a gap below $\k$, then $\k$ was $\l$-supercompact in the ground model, and further, every supercompactness measure in the extension extends a supercompactness measure in the ground model.

\proof If $j:V[G]\to M[j(G)]$ is the ultrapower by the $\l$-supercompactness measure $\mu$ then the Gap Forcing Theorem implies that the restricted embedding $j:V\to M$ is definable from parameters in $V$, and Lemma \Assist\ implies that $M^\l\of M$ in $V$. In particular, $[\id]_\mu=j\image\l\in M$, and so $\mu\intersect V$ must be in $V$, as desired.\qed  

\corollary. Gap forcing creates no new almost huge cardinals, huge cardinals, or $n$-huge cardinals for any $n\in\omega$. 

\proof This argument is just the same. If $j:V[G]\to M[j(G)]$ is an almost hugeness or hugeness embedding in $V[G]$, then the Gap Forcing Theorem implies that the restricted embedding $j:V\to M$ is definable from parameters in $V$ and Lemma \Assist\ shows that it is has the corresponding amount of hugeness there.\qed

Let me close with the following observation.

\theorem Observation. The closure assumption on the embedding in the Gap Forcing Theorem cannot be omitted, because if there are two normal measures on the measurable cardinal $\k$ in $V$ then after merely adding a Cohen real $x$ there is an embedding $j:V[x]\to M[x]$ that does not lift an embedding from the ground model.

\proof Suppose that $\mu_0$ and $\mu_1$ are normal measures on $\k$ in $V$ and $x$ is a $V$-generic Cohen real. By the Levy-Solovay Theorem \cite[LevSol67], these measures extend uniquely to measures $\bar\mu_0$ and $\bar\mu_1$ in $V[x]$, and furthermore the ultrapowers by the measures $\bar\mu_0$ and $\bar\mu_1$ in $V[x]$ are the unique lifts of the corresponding ultrapowers by $\mu_0$ and $\mu_1$ in $V$. Let $j:V[x]\to M[x]$ be the $\w$-iteration determined in $V[x]$ by selecting at the $n^\th$ step either the the image of $\bar\mu_0$ or of $\bar\mu_1$, respectively, depending on the $n^\th$ digit of $x$. If $\<\k_n\st n<\w>$ is the critical sequence of this embedding, then for any $X\of\k$ the standard arguments show that $\k_n\in j(X)$ if and only if $X$ is in the measure whose image is used at the $n^\th$ step of the iteration. Suppose now towards a contradiction that the restricted embedding $j\restrict V$ is amenable to $V$. I will show that from $j\restrict P(\k)^V$ one can iteratively recover the digits of $x$. First, by computing in $V$ the set $\set{X\of\k\st \k\in j(X)}$, we learn which measure was used at the initial step of the iteration and thereby also learn the initial digit of $x$. This information also tells us the value of $\k_1=j_{\mu_{x(0)}}(\k)$. Continuing, we can compute in $V$ the set $\set{X\of\k\st \k_1\in j(X)}$ to know the next measure that was used and thereby learn the next digit of $x$ and the value of $\k_2$, and so on. Thus, from $j\restrict P(\k)^V$ in $V$ we would be able to recursively recover $x$, contradicing the fact that $x$ is not in $V$.\qed 

The argument works equally well with any small forcing; one simply uses a longer iteration. 

\medskip
{\parindent=0pt\tenpoint\tightlineskip\sc 
Mathematics, City University of New York, College of Staten Island\par 
2800 Victory Boulevard, Staten Island, NY 10314\par
\tt hamkins@math.csi.cuny.edu\par
http://www.math.csi.cuny.edu/$\sim$hamkins\par}

\quiet\section Bibliography

\nopagenumbers
\parindent=0pt
\newbox\Article
\newbox\Journal
\newbox\Author
\newbox\Vol
\newbox\No
\newbox\Year
\newbox\Page
\newbox\Book
\newbox\Publisher
\newbox\Pubaddr
\newbox\Key
\newbox\Editor
\newbox\Comment
\newbox\Note
\def\entry#1#2\par{\item{#1\quad}\hskip-1.1em#2\par}
\def\article#1{\setbox\Article=\hbox{\sl #1, }}
\def\journal#1{\setbox\Journal=\hbox{\rm #1 }}
\def\author#1{\setbox\Author=\hbox{\sc #1, }}
\def\vol#1{\setbox\Vol=\hbox{\bf #1 }}
\def\no#1{\setbox\No=\hbox{no. #1 }}
\def\year#1{\setbox\Year=\hbox{\rm({\oldstyle #1}) }}
\def\page#1{\setbox\Page=\hbox{\rm p. #1 }}
\def\book#1{\setbox\Book=\hbox{\it #1, }}
\def\publisher#1{\setbox\Publisher=\hbox{\rm #1, }}
\def\pubaddr#1{\setbox\Pubaddr=\hbox{\rm #1, }}
\def\key#1{\setbox\Key=\hbox{#1}}
\def\editor#1{\setbox\Editor=\hbox{\rm(#1, Ed.) }}
\def\comment#1{\setbox\Comment=\hbox{\rm #1}}
\def\note#1{\setbox\Note=\hbox{\rm #1 }}
\def\ref#1\par{\smallskip{#1
  \entry{\ifhbox\Key\unhbox\Key\else[\ ]\fi}%
  \unhbox\Author\unhbox\Note
  \ifhbox\Book \unhbox\Book\unhbox\Publisher\unhbox\Pubaddr
               \unhbox\Editor\unhbox\Page\unhbox\Year\unhbox\Comment
  \else \unhbox\Article\unhbox\Journal\unhbox\Vol\unhbox\No\unhbox\Editor
        \unhbox\Page\unhbox\Year\unhbox\Comment\fi\par}}

\tenpoint\tightlineskip

\ref
\author{Joel David Hamkins}
\article{Canonical seeds and Prikry trees}
\journal{Journal of Symbolic Logic}
\vol{62}
\no{2}
\page{373-396}
\year{1997}
\key{[Ham97]}

\ref
\author{Joel David Hamkins}
\article{Destruction or preservation as you like it}
\journal{Annals of Pure and Applied Logic}
\year{1998}
\vol{91}
\page{191-229}
\key{[Ham98a]}

\ref
\author{Joel David Hamkins}
\article{Small forcing makes any cardinal superdestructible}
\journal{Journal of Symbolic Logic}
\year{1998}
\vol{63}
\no{1}
\page{51-58}
\key{[Ham98b]}

\ref
\author{Joel David Hamkins and Saharon Shelah}
\article{Superdestructibility: a dual to the Laver preparation}
\journal{Journal of Symbolic Logic}
\year{1998}
\vol{63}
\no{2}
\page{549-554}
\key{[HamShl98]}

\ref
\author{Joel David Hamkins}
\article{Gap forcing: generalizing the Levy-Solovay theorem}
\journal{submitted to the Bulletin of Symbolic Logic}
\key{[Ham$\infty$]}

\ref
\author{Joel David Hamkins and W. Hugh Woodin}
\article{Small forcing creates neither strong nor Woodin cardinals}
\comment{forthcoming}
\key{[HamWdn]}

\ref
\author{Richard Laver}
\article{Making the supercompactness of $\kappa$ indestructible under 
 $\kappa$-directed closed forcing}
\journal{Israel Journal Math}
\vol{29}
\year{1978}
\page{385-388}
\key{[Lav78]}

\ref
\author{Dana S. Scott}
\article{Measurable cardinals and constructible sets}
\journal{Bulletin of the Polish Academy of Sciences, Mathematics}
\vol{9}
\year{1961}
\page{521-524}
\key{[Sco61]}

\ref
\author{Jack Silver}
\article{The Consistency of the Generalized Continuum Hypothesis with the existence of a Measurable Cardinal}
\journal{Axiomatica Set Theory, Proc. Symp. Pure Math. 13 American Mathematical Society}
\vol{I}
\editor{D. Scott}
\year{1971}
\page{383-390}
\key{[Sil71]}

\ref
\author{Levy Solovay}
\article{Measurable cardinals and the Continuum Hypothesis}
\journal{IJM}
\vol{5}
\year{1967}
\page{234-248}
\key{[LevSol67]}

\bye